\documentclass[letterpaper, 10 pt, conference]{ieeeconf}
\IEEEoverridecommandlockouts

\usepackage{graphicx}
\usepackage{amsmath}

\usepackage{cite}
\usepackage{authblk}
\usepackage{amssymb}
\usepackage{subfigure}
\usepackage{multirow}
\usepackage{color}
\usepackage{tikz} \usetikzlibrary{arrows,calc,decorations.pathreplacing}
         
\newtheorem{theorem}{Theorem}

\newtheorem{corollary}{Corollary}

\title{A Differential Game Approach for Beyond Visual Range Tactics}
\author{Eloy Garcia, David W. Casbeer, Dzung Tran, and Meir Pachter 
\thanks{This work has been supported in part by AFOSR LRIR No. 18RQCOR036.}
\thanks{E. Garcia, D. Casbeer, and D. Tran are with the Control Science Center of Excellence, Air Force Research Laboratory, Wright-Patterson AFB, OH 45433.  Corresponding author \ttfamily{eloy.garcia.2@us.af.mil}}
\thanks{M. Pachter is with the Department of Electrical Engineering, Air Force Institute of Technology, Wright-Patterson AFB, OH 45433.}}

\begin{document}
\maketitle 
                                         
\begin{abstract}
An operational relevant conflict between teams of autonomous vehicles in the Beyond Visual Range domain is addressed in this paper. 
Optimal strategies are designed in order for a team of air interceptors to protect a high value asset and block the attacking team at a safe distance from such asset. The attacking agents take specific roles of leader and wingman and also devise their own optimal strategies in order to launch an attack as close as possible from the asset. The problem is formulated as a zero-sum differential game between players with different speed over two stages: the attack and the retreat stages. For each stage the state-feedback optimal strategies of each player are derived in analytical form.
\end{abstract}

\section{Introduction} \label{sec:intro}

In the military, a migration is underway from conventional forces into unmanned, modular, cooperative, and consequently more manageable teams of assets. This process requires additional planning and incorporation of new models, methods, and team strategies that leverage the new capabilities and properties of these autonomous systems \cite{chandler2004cooperative,Archibald08,wang2017cooperative,Beard02,Castanon09,Lechevin09,Garcia16,Ernest15}. 
It is also expected that a team of intelligent and autonomous assets will face an opposing intelligent team with similar capabilities and the same flexibility to leverage cooperation on their own. This combat scenario will require efficient command, control, and guidance strategies which allow a team of autonomous systems to prevail and complete specified missions in the presence of adversarial units. Strategies which enable coordination and cooperation among teammates should be devised as well in order to determine the best course of action.

Mathematical analysis of conflict is performed by leveraging game theory. Game theory techniques are used to design optimal strategies that are also robust to unknown adversary's actions. Controllers and decision making algorithms designed based on this approach are able to adapt to potential enemy actions in a rapidly changing battle space and to take advantage of enemy's deviation from optimal solutions.  
The authors of \cite{galati2008near} emphasized the need for strategies based on game theoretic analysis to address automated battle scenarios.  The key aspect is to synthesize Nash-equilibrium strategies where players do not obtain any improvement in their performance by deviating from these strategies. The papers  \cite{liu2003application,galati2007effectiveness,ganapathy2003agreement} offer similar game theory techniques for battle management. However, the previous references concerning game theory analysis of combat operations employed static game formulations. In order to formally address the dynamical aspect exhibited in autonomous air battle management, the tools and methods of differential game theory are needed  \cite{Isaacs65,weintraub20}. Differential game theory is concerned with the analysis of non-cooperative conflicts where the underlying processes are governed by differential equations. Differential game theory has been successfully applied for motion planning \cite{zhou2018efficient,chen2014path} and, especially, for analysis of pursuit-evasion games \cite{Coon17,Bopardikar09,Garcia2019,Fuchs17,Pachter17,fisac2015pursuit,Weintraub20IET,GarciaCDC18}.

The problem considered is an attack-defense of high-value stationary asset where a team of two Unmanned Aerial Vehicles (UAV) is tasked to attack it while a team of two interceptors tries to deny access and protect the valuable asset.
The problem takes place in the Beyond Visual Range (BVR) domain. This is a highly relevant and practical scenario involving several agents.
We denote the attacking UAVs as the blue team and the interceptors as the red team. 
The problem is addressed in two stages, the attack and the retreat stages. We pay particular attention to safety and survivability of the blue team; hence, as they choose the mode of attack, the leader-wingman roles are assigned. In the attack stage, the leader's role is to penetrate the area protected by the red interceptors and attack the red entities while the wingman flies in formation. In the retreat stage, the wingman is in charge of protecting the leader while the latter is trying to evade the weapons fired by the red interceptors. It is important to note that due to this type of cooperation between blue vehicles, it is possible for the leader to penetrate deeper into the area protected by red interceptors knowing that its teammate or wingman will aid him to escape. Consequently, traditional weapon engagement zones (WEZ) become less relevant in the design of cooperative tactics.

This problem was originally proposed in \cite{Garcia21GNC} where a preliminary analysis was presented. The analysis in \cite{Garcia21GNC} considers all air vehicles, blue UAVs and red interceptors, to have the same speed. Also, the retreat stage was not formally analyzed in that reference. In this paper we generalize those preliminary results in order to consider players with different speed. This generalization has important practical applications by addressing the players' sense of urgency in the way human pilots usually approach such concerns. In addition, in this paper, we provide a detailed analysis of the retreat stage. Active defense of an aircraft is analyzed in the presence of two attacker missiles and two defender missiles. Finally, the developed autonomous battle management method is exemplified in simulations. The sense of urgency in terms of survivability of the blue team is incorporated and the optimal strategies, which are obtained based on simple motion models, are implemented in aircraft models that account for both turning rate and acceleration constraints. 

The rest of the paper is organized as follows. The conflict between the blue and the red teams in the BVR domain is formulated as a differential game in Section \ref{sec:PD}. Section \ref{sec:Coop} analyzes the attack stage where the air vehicles have different speed and cooperative optimal strategies for attack-defense of the high value asset are developed.
In Section \ref{sec:Retreat} optimal strategies for the retreat stage are presented.
Section \ref{sec:Ex} provides illustrative examples and conclusions follow in Section \ref{sec:concl}.

 \section{The BVR Conflict} \label{sec:PD}
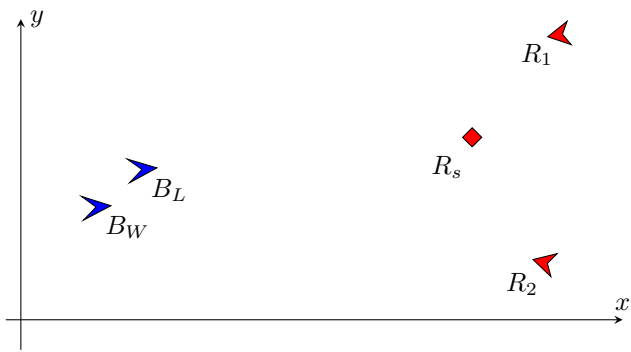
\begin{figure}[htb]\centering
\begin{tikzpicture}[>=stealth]
	\tikzset{
	   uav/.pic = {  
	      \draw [pic actions] (0,0)
	         -- (-1mm,1.75mm)
	         -- ++(3.2mm,-1.75mm)  
	         -- ++(-3.2mm,-1.75mm)
	         -- (0,0);
	   }
         }
	\tikzset{
	   missile/.pic = {  
	      \draw [pic actions] (0,0)
	         -- (-2mm,1.75mm)
	         -- ++(4.2mm,-1.75mm)  
	         -- ++(-4.2mm,-1.75mm)
	         -- (0,0);
	   }
	}
	\tikzset{
	   asset/.pic = {  
	      \draw [pic actions] (0,0)
	         -- (2mm,0mm)
	         -- ++(0mm,2mm)  
	         -- ++(-2mm,0mm)
	         -- (0,0);
	   }
	}
	\coordinate (R2) at (7,.75);  
	\coordinate (R1) at (7.2,3.8);
	\coordinate (BL) at (1.61,2);
	\coordinate (BW) at (1.0,1.5);
	\coordinate (A) at (6,2.3);
	\draw (R2) pic [rotate=165,fill=red,scale=.9] {uav};  
	\draw (R2) node [below left] {$R_2$};
	\draw (BL) pic [rotate=6,fill=blue,scale=.9] {missile};
	\draw (BL) node [below right] {$B_L$};
	\draw (R1) pic [rotate=-170,fill=red,scale=.9] {uav};
	\draw (R1) node [below left] {$R_1$};
	\draw (BW) pic [rotate=5,fill=blue,scale=.9] {missile};
	\draw (BW) node [below right] {$B_W$};
	\draw (A) pic [rotate=45,fill=red,scale=.9] {asset};
	\draw (A) node [below left] {$R_s$};
	\node [right] (Y) at (0,4) {$y$};
	\draw [->](0,-.4) -- (0,4);
	\node [above] (X) at (8,0) {$x$};
	\draw [->](-.2,0) -- (8,0);
 \end{tikzpicture}
\caption{BVRT: Attack Stage}
\label{fig:AppSta}
\end{figure}

A conflict between two UAVs, which are labeled as $B_L$ and $B_W$, and two air interceptors, which are labeled as $R_1$ and $R_2$, is considered. The stationary asset is denoted by $R_s$. The blue team's ($B_L$ and $B_W$) objective is to penetrate the area protected by the red team ($R_1$ and $R_2$) and reach within firing distance, $\rho_s>0$, of $R_s$. The red interceptors are tasked to block the blue team of UAVs by reaching a firing range $\rho>0$ with respect to a blue UAV. A red interceptor will fire an air-to-air missile against a blue UAV located within firing range $\rho$. The blue UAVs place a high value to their own safety and, similar to the red interceptors, are equipped with air-to-air missiles of similar range $\rho$. Hence, if a blue UAV is within range of a red interceptor it will fire a missile and then retreat to evade the red interceptor's attack.  This engagement is illustrated in Fig. \ref{fig:AppSta}. The missiles' range $\rho$ is typically longer than visual range, hence, the problem is to devise corresponding tactics for a BVR combat scenario.

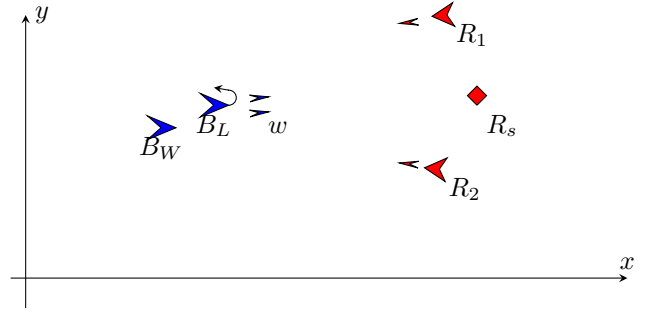
\begin{figure}[htb]\centering
\begin{tikzpicture}[>=stealth]
	\tikzset{
	   uav/.pic = {  
	      \draw [pic actions] (0,0)
	         -- (-1mm,1.75mm)
	         -- ++(3.2mm,-1.75mm)  
	         -- ++(-3.2mm,-1.75mm)
	         -- (0,0);
	   }
         }
	\tikzset{
	   missile/.pic = {  
	      \draw [pic actions] (0,0)
	         -- (-2mm,1.75mm)
	         -- ++(4.2mm,-1.75mm)  
	         -- ++(-4.2mm,-1.75mm)
	         -- (0,0);
	   }
	}
	\tikzset{
	   missile2/.pic = {  
	      \draw [pic actions] (0,0)
	         -- (-2mm,.75mm)
	         -- ++(4.2mm,-.75mm)  
	         -- ++(-4.2mm,-.75mm)
	         -- (0,0);
	   }
	}
	\tikzset{
	   asset/.pic = {  
	      \draw [pic actions] (0,0)
	         -- (2mm,0mm)
	         -- ++(0mm,2mm)  
	         -- ++(-2mm,0mm)
	         -- (0,0);
	   }
	}
	\coordinate (R2) at (5.5,1.45);  
	\coordinate (R1) at (5.6,3.5);
	\coordinate (BL) at (2.5,2.3);
	\coordinate (BW) at (1.8,2.0);
	\coordinate (A) at (6,2.3);
	\coordinate (M1) at (3.1,2.4);
	\coordinate (M2) at (3.1,2.2);
	\coordinate (M3) at (5.1,3.4);
	\coordinate (M4) at (5.1,1.52);
	\draw (R2) pic [rotate=175,fill=red,scale=.9] {uav};  
	\draw (R2) node [below right] {$R_2$};
	\draw [->] 
			(BL) ++(0:.2) arc [start angle=-90, end angle=95, radius=.1cm]
	 	       (2.7,2.5) -- ++(170:.2);
	\draw (BL) pic [rotate=0,fill=blue,scale=.9] {missile};
	\draw (BL) node [below] {$B_L$};
	\draw (R1) pic [rotate=185,fill=red,scale=.9] {uav};
	\draw (R1) node [below right] {$R_1$};
	\draw (BW) pic [rotate=0,fill=blue,scale=.9] {missile};
	\draw (BW) node [below] {$B_W$};
	\draw (A) pic [rotate=45,fill=red,scale=.9] {asset};
	\draw (A) node [below right] {$R_s$};
	\draw (M1) pic [rotate=5,fill=blue,scale=.6] {missile2};
	\draw (M2) pic [rotate=5,fill=blue,scale=.6] {missile2};
	\draw (M2) node [below right] {$w$};
	\draw (M3) pic [rotate=185,fill=red,scale=.6] {missile2};
	\draw (M4) pic [rotate=175,fill=red,scale=.6] {missile2};
	\node [right] (Y) at (0,3.5) {$y$};
	\draw [->](0,-.4) -- (0,3.5);
	\node [above] (X) at (8,0) {$x$};
	\draw [->](-.2,0) -- (8,0);
 \end{tikzpicture}
\caption{BVRT: Retreat Stage}
\label{fig:ReSta}
\end{figure}

The blue team strives to reach as close as possible to the red asset $R_s$ in the presence of intelligent adversaries, namely, the red interceptors. The red interceptors, in turn, try to protect the asset and block the blue team as far as possible from $R_s$. In order to balance lethality and safety, the scenario consists of two stages: the attack stage and the retreat stage. 
In the attack stage, $B_L$ engages either $R_S$ or the red interceptors. $B_L$ immediately retreats after firing weapons, whether it is against the interceptors or $R_S$. This sequence is shown in Fig. \ref{fig:ReSta}. In the retreat stage $B_W$ takes the role of protecting $B_L$ by firing defender missiles to intercept the missiles fired by the red interceptors against $B_L$. 


In general, we consider $\rho_s\neq \rho$ to allow the case where $B_L$ is equipped with two different types of weapons: one to launch against aircraft interceptors and one to launch against static ground assets.
In this work we assume an engagement range $\rho$ instead of a weapon engagement zone (WEZ) and we also assume simple motion models for the blue and red aircraft. Both assumptions are due to the nature of the BVR problem under consideration as opposed to a close-range conflict where dynamic constraints play a more significant role. More importantly, by employing cooperative active target defense strategies in the retreat stage, it is possible for $B_L$ to penetrate deeper into the area protected by red interceptors knowing that its teammate or wingman will aid him to escape; thus, the WEZ becomes less relevant. This means that $\rho$ could be selected to be smaller than air-to-air missile range and in accordance to active target defense strategies. This level of cooperation between $B_L$  and $B_W$ allows for tactics and strategies that are more aggressive and lethal while avoiding loss of blue aircraft. Additionally, in Section \ref{sec:Ex}, the state-feedback strategies obtained based on simple kinematics will be used in aircraft models with both turning rate and acceleration constraints. Those examples highlight the great importance of synthesizing state-feedback strategies for dynamic conflicts and the applicability and effectiveness of these strategies in BVR scenarios.

\section{The Attack Stage} \label{sec:Coop}

\subsection{Problem Formulation}   \label{subsec:BVRTdg}
\begin{figure}
	\begin{center}
		\includegraphics[width=8cm,trim=.9cm .5cm 1.5cm .2cm]{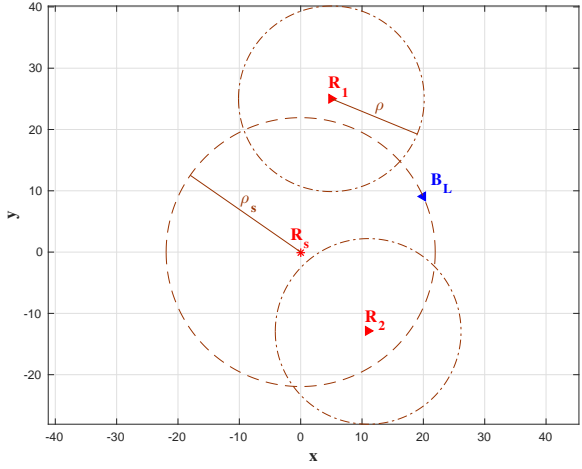}
	\caption{$B_L$ reaches a distance $\rho_s$ with respect to $R_s$ before being blocked by any interceptor}
	\label{fig:rhoA}
	\end{center}
\end{figure}

In this section the attack stage of the problem is formulated as a zero-sum differential game between the blue and the red teams. 
In the attack stage the players are Blue Lead $B_L$ and two red interceptors $R_1$ and $R_2$ who protect the stationary asset $R_s$. $B_L$ strives to reach an engagement distance $\rho_s>0$ with respect to $R_s$ and fire a missile/weapon to destroy $R_s$. The interceptors $R_1$ and $R_2$ are tasked to block $B_L$ and deny penetration of $B_L$ within a distance $\rho_s$ of $R_s$. If either one of $R_1$ or $R_2$ is able to achieve an engagement range $\rho>0$ with respect to $B_L$ they will block $B_L$ by firing a missile aimed at $B_L$. If $B_L$ is blocked by $R_1$ or $R_2$, he will also launch a missile against the interceptor located within range $\rho$ and retreat immediately since the interceptor has launched a missile at him. On the other hand, if $B_L$ is able to penetrate and reach a distance $\rho_s$ from $R_s$ then he will fire a weapon aimed at $R_s$ followed by retreat. For instance, Fig. \ref{fig:rhoA} shows an example where $B_L$ is able to reach a distance $\rho_s$ from $R_s$ before being blocked by any interceptor. 

Since $R_s=(x_s,y_s)$ is static, the engagement zone of $R_s$ is defined as follows
\begin{align}
	\Gamma =\{ x,y | \sqrt{(x-x_s)^2+(y-y_s)^2}\leq \rho_s\}.
\end{align}
The states of $B_L$, $R_1$, and $R_2$ are respectively specified by their Cartesian coordinates $\textbf{x}_{B}=(x_B,y_B)$, $\textbf{x}_1=(x_1,y_1)$ and $\textbf{x}_2=(x_2,y_2)$. The players have constant speeds $v_B$, $v_1$, and $v_2$.
The complete state of the differential game is defined by $\textbf{x}:=(x_B, y_B, x_1, y_1, x_2, y_2)\in \mathbb{R}^6$. 
The  control variable of $B_L$ is its instantaneous heading angle, $\textbf{u}_B=\left\{\theta_B\right\}$. The interceptors affect the state of the game by choosing their instantaneous headings, $\textbf{u}_R=\left\{\theta_1,\theta_2\right\}$. The dynamics $\dot{\textbf{x}}=\textbf{f}(\textbf{x},\textbf{u}_B,\textbf{u}_R)$ are specified by the system of linear differential equations
\begin{align}
 \left.
	\begin{array}{l l}
        \dot{x}_B&=v_B\cos\theta_B,    \ \ \ \ \ x_B(0)=x_{B_0}   \\
	\dot{y}_B&=v_B\sin\theta_B,  \ \ \ \ \ y_B(0)=y_{B_0}  \\
	\dot{x}_1&=v_1\cos\theta_1,   \ \ \ \ \ \ x_1(0)=x_{1_0}  \\
	\dot{y}_1&=v_1\sin\theta_1,   \ \ \ \ \ \ y_1(0)=y_{1_0}  \\
	\dot{x}_2&=v_2\cos\theta_2,   \ \ \ \ \ \ x_2(0)=x_{2_0}  \\
	\dot{y}_2&=v_2\sin\theta_2,   \ \ \ \ \ \ y_2(0)=y_{2_0}
	\end{array}  \right.   \label{eq:xT}
\end{align}
where the admissible controls are the players' headings $\theta_B, \theta_1, \theta_2 \in [-\pi,\pi)$. 
The initial state of the system is defined as
\begin{align}
	\textbf{x}_0 := ( x_{B_0}, y_{B_0}, x_{1_0}, y_{1_0}, x_{2_0}, y_{2_0}) = \textbf{x}(t_0).  \nonumber
\end{align}
In the attack stage we assume that $B_L$ flies at low speed $v_B<v_1=v_2$. Define the speed ratio $\beta=v_1/v_B>1$. Without loss of generality, the speeds of the players are normalized so $v_1=v_2=\beta$ and $v_B=1$. 

In this paper we consider the case where $B_L$ implements a different speed than the red interceptors. Besides considering players which do not necessarily have the same speed, this work allows to implement tactics that consider survivability and a sense of urgency as a key aspect. By flying at low speed in the attack stage, $B_L$ is able to increase maneuverability at the start of the retreat stage, then, $B_L$ increases its speed, and successfully escapes from the attack by the red interceptors.

The target (or termination) set of the attack stage is 
\begin{align}
   \mathcal{T} :=  \mathcal{T}_\Gamma    \   \bigcup \   \mathcal{T}_R   \label{eq:TwoSets}
\end{align}
where 
\begin{align}
   \mathcal{T}_\Gamma:=   \big\{ \ \textbf{x} \ |  \sqrt{ (x_B-x_s)^2+(y_B-y_s)^2} = \rho_s \big\}   \label{eq:Set1}
\end{align}
represents the outcome where $B_L$ is able to reach the engagement zone $\Gamma$ before being blocked by the interceptors. On the other hand
\begin{align}
   \left.
	\begin{array}{l l}
	 \mathcal{T}_R:= \big\{ \ \textbf{x} \ | \sqrt{ (x_B-x_1)^2 + (y_B-y_1)^2} = \rho \big\}  \\
   \qquad   \quad  \bigcup \    \big\{ \ \textbf{x} \ | \sqrt{ (x_B-x_2)^2 + (y_B-y_2)^2} = \rho \big\} \label{eq:Set2}
\end{array}  \right. 
\end{align}
represents the outcome where  at least one of the interceptors is able to block $B_L$ from penetrating the defended area and engage the asset $R_s$. 

The differential game with termination set as given in \eqref{eq:TwoSets} belongs to the class of two termination set differential games \cite{getz1979qualitative,Getz1981capturability}. The two termination set differential game concept was introduced in order to extend classical pursuit-evasion games where only one termination set is contemplated.
For instance, the pursuer tries to minimize the cost to reach the termination set whereas the evader wants to maximize the payoff to reach the termination set or, when possible, to avoid reaching that set at all. 
The two termination set differential game is useful in the analysis of combat games \cite{ardema1985combat} where the roles of pursuer and evader are not designated ahead of time; instead each player wants to defeat the opponent by terminating the game in its own termination set. In this problem, $B_L$ strives to reach the engagement zone $\Gamma$ before being within distance $\rho$ from any of the interceptors.

An important feature of this work is to leverage cooperation between $R_1$ and $R_2$. In certain cases, it is possible that individual optimal solutions of each interceptor against $B_L$ do not block $B_L$ from reaching the engagement zone $\Gamma$. However, if the interceptors cooperate they can find an strategy that will successfully block $B_L$ from entering $\Gamma$.
It will be shown that, depending on the initial conditions, the optimal strategy is for $R_1$ and $R_2$ to cooperate and block $B_L$ simultaneously. In other words, the interceptors can block $B_L$ the farthest from $R_s$ when they achieve the engagement range $\rho$ exactly at the same time instant. The optimal strategy of $B_L$ in such a case is to follow the same solution in order to reach as close to $R_s$ as possible. In the case where $B_L$ reaches the engagement range $\rho$ simultaneously with respect to $R_1$ and $R_2$, he will fire a weapon against each interceptor and retreat immediately.


Let $\mathcal{R}_R\in \mathbb{R}^6$ denote the red team winning subspace where, under optimal play, the terminal condition is \eqref{eq:Set2}.
The terminal time $t_f$ is defined as the time instant when the state of the system satisfies \eqref{eq:Set2}, at which time the terminal state is $\textbf{x}_f: = ( x_{B_f}, y_{B_f}, x_{1_f}, y_{1_f}, x_{2_f}, y_{2_f}) = \textbf{x}(t_f)$.
The terminal cost/payoff functional is
\begin{align}
  J(\textbf{u}_B(t),\textbf{u}_R(t);\textbf{x}_0)=\Phi(\textbf{x}_f)	\label{eq:costDG2}
\end{align}
where
\begin{align}
  \Phi(\textbf{x}_f):=\sqrt{ (x_{B_f}-x_s)^2+(y_{B_f}-y_s)^2}. 	\label{eq:costDG}
\end{align}
The cost/payoff functional depends only on the terminal state - the capture game is a terminal cost/Mayer type game. Its Value is given by
\begin{align}
  V(\textbf{x}_0):= \min_{\textbf{u}_B(\cdot)} \ \max_{\textbf{u}_R(\cdot)} J(\textbf{u}_B(\cdot),\textbf{u}_R(\cdot);\textbf{x}_0)	\label{eq:costDG3}
\end{align}
subject to \eqref{eq:xT} and \eqref{eq:Set2}, where $\textbf{u}_B(\cdot)$ and $\textbf{u}_R(\cdot)$ are the players' state feedback strategies.

The performance functional \eqref{eq:costDG} is an important measure of risk associated to this combat scenario. The interceptors not only want to block $B_L$ from reaching $\Gamma$ but they also want to maximize the distance between $B_L$ terminal position and the location of the protected asset $R_s$. Naturally, $B_L$ aims at minimizing the same distance in order to approach as close as possible to $R_s$ before being forced to retreat. The key feature of choosing the performance functional \eqref{eq:costDG} is that under non-optimal play by one of the teams, the adversary can see an improvement in its performance and actually change the outcome of the game.  The state-feedback saddle-point strategies of this differential game are necessary in order to obtain the Nash equilibrium properties. Then,  $B_L$ is guaranteed to further decrease its terminal distance with respect to $R_s$ if the interceptors fail to implement their optimal strategy. Under non-optimal play by the interceptors, $B_L$ could potentially be able to reach $\Gamma$ and actually win the game by engaging the asset $R_s$. On the contrary, if $B_L$ does not implement its optimal strategy, the interceptors will be able to block him farther away from $R_s$ reducing, in this way, the threat to the protected asset $R_s$.

\subsection{Cooperative Tactical Operations in the Attack Stage} \label{sec:CoopStr}

In this section we consider the case where two red interceptors $R_1$ and $R_2$ wish to simultaneously combine forces and block $B_L$ as far as possible from $R_s$. 

\begin{theorem}  \label{th:SL}
Consider the BVR differential game and assume that $\textbf{x}\in\mathcal{R}_R$. The headings of the players $B_L$, $R_1$ and $R_2$ are constant under optimal play and the optimal trajectories are straight lines.
\end{theorem}
\textit{Proof}. The proof is provided in Appendix A.

\begin{figure}
	\begin{center}
		\includegraphics[width=8cm,trim=.9cm .4cm 1.5cm .5cm]{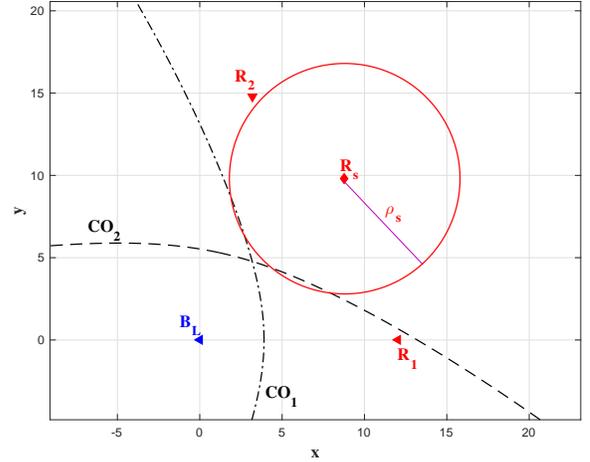}
	\caption{Relative coordinate frame where $B_L=(0,0)$, $R_1=(x'_1,0)$, $R_2=(x'_2,y'_2)$, and $R_s=(x'_s,y'_s)$.}
	\label{fig:TwoCos}
	\end{center}
\end{figure}

Since the optimal headings of each player are constant, the dominance region of $B_L$ with respect to each interceptor is delineated by a Cartesian Oval (CO). We consider, without loss of generality, the relative coordinate frame shown in Fig. \ref{fig:TwoCos}, where $B_L$ is located at the origin.  In this coordinate frame we have that $R_1$ is located on the positive $x$-axis; its coordinates are given by $R_1=(x'_1,0)$ where 
 \begin{align}
\left.
	\begin{array}{l l}
  	x'_1=\sqrt{(x_1-x_B)^2+(y_1-y_B)^2} \\
         \end{array}   \right.  \nonumber 
\end{align}
The positions of $R_2$ and $R_s$ are given by $R_2=(x'_2,y'_2)$ and $R_s=(x'_s,y'_s)$ where
 \begin{align}
\left.
	\begin{array}{l l}
	x'_2=d_2\cos(\lambda_2-\lambda_1) \\
	y'_2=d_2\sin(\lambda_2-\lambda_1)  \\
	x'_s=d_s\cos(\lambda_s-\lambda_1) \\
	y'_s=d_s\sin(\lambda_s-\lambda_1) 
         \end{array}   \right.  \nonumber 
\end{align}
$d_2=\sqrt{(x_2-x_B)^2+(y_2-y_B)^2}$, $\lambda_2=\arctan(\frac{y_2-y_B}{x_2-x_B})$, $d_s=\sqrt{(x_s-x_B)^2+(y_s-y_B)^2}$, $\lambda_s=\arctan(\frac{y_s-y_B}{x_s-x_B})$, and $\lambda_1=\arctan(\frac{y_1-y_B}{x_1-x_B})$. In the relative coordinate frame, the CO between $B_L$ and $R_1$ is given by  
\begin{align}
	 \sqrt{(x-x'_1)^2+y^2}=\rho + \beta \sqrt{x^2+y^2}   \label{eq:CO1}
\end{align}
where $\beta>1$ is the speed ratio and $\rho>0$ is the engagement distance between aircraft. 
Similarly, the CO between $B_L$ and $R_2$ is given by  
\begin{align}
	 \sqrt{(x-x'_2)^2+(y-y'_2)^2}=\rho + \beta \sqrt{x^2+y^2}   \label{eq:CO2}
\end{align}
In the presence of two interceptors, the dominance region of $B_L$ is given by the intersection of the two COs \eqref{eq:CO1} and \eqref{eq:CO2}.  Let $\mathcal{B}_L$ denote the dominance region of $B_L$. 

Figure \ref{fig:TwoCos} shows an interesting case where cooperation between red interceptors is necessary in order to successfully block $B_L$. For instance, if only one of the interceptors commits to block $B_L$ there is an strategy for $B_L$ to reach $\Gamma$ and being able to engage $R_s$. This is due to the fact that each CO intersects the engagement zone $\Gamma$; however, if both interceptors cooperate in order to block $B_L$, then, they are able to constrain $B_L$'s dominance region. In this example, the closest that $B_L$ can reach with respect to $\Gamma$ is the intersection point of both COs.

Let us define the orthogonal bisector of the segment $\overline{R_1R_2}$ which is given by the equation $y=mx+n$, where
\begin{align}
 \left.
	 \begin{array}{l l}
	m=\frac{x'_1-x'_2}{y'_2}, \ \ \ \ \
	n= -\frac{1}{2} \frac{x_1^{'2}-y_2^{'2}-x_2^{'2}}{y'_2}.
   \label{eq:mni}
\end{array}  \right.  
\end{align}
In the following theorem we assume that simultaneous capture is the optimal play; then, the line $y=mx+n$ passes through the intersection of the two COs \eqref{eq:CO1} and \eqref{eq:CO2}. Define the functions
\begin{align}
 \left.
	 \begin{array}{l l}
	O(x) =\sqrt{(x-x'_1)^2+(mx+n)^2} \\
	\qquad \quad - \beta \sqrt{x^2+(mx+n)^2} -\rho
   \label{eq:OvalFn}
\end{array}  \right.  
\end{align}
and
\begin{align}
 \left.
	 \begin{array}{l l}
	J(x) = \sqrt{(x'_s-x)^2+(y'_s-mx-n)^2}
   \label{eq:CostPair}
\end{array}  \right.  
\end{align}

\begin{theorem}  \label{th:coop}
\textbf{Cooperative blocking strategy}. Consider the attack stage of the BVR differential game and assume that $\textbf{x}\in \mathcal{R}_{R_s}$. 
 The optimal headings of players $B_L$, $R_1$, and $R_2$ in the relative coordinate frame are given by 
 \begin{align}
 \left.
	 \begin{array}{l l}
	  \cos\theta_B^* &= \frac{x^*}{\sqrt{x^{*2} + y^{*2}}}      \\
	    \sin\theta_B^* &= \frac{y^*}{\sqrt{x^{*2} + y^{*2}}}      \\
	     \cos\theta_1^* &= \frac{x^*-x'_1}{\sqrt{(x^*-x'_1)^2 + y^{*2}}}      \\
	          \sin\theta_1^* &= \frac{y^*}{\sqrt{(x^*-x'_1)^2 + y^{*2}}}   \\
	           \cos\theta_2^* &= \frac{x^*-x'_2}{\sqrt{(x^*-x'_2)^2 + (y^*-y'_2)^2}}      \\
	          \sin\theta_2^* &= \frac{y^*-y'_2}{\sqrt{(x^*-x'_2)^2 + (y^*-y'_2)^2}} 
\end{array}  \right.    \label{eq:OptimalInputsPCRcoop}
\end{align}
 where $y^*=mx^*+n$ and $x^*$ is given by
 \begin{align}
 \left.
	 \begin{array}{l l}
	x^*=\arg\min \{ J(x_{in_1}),J(x_{in_2})\}
	\end{array}  \right.   \label{eq:OptXpair}
\end{align}
 where $x_{in_1}$ and $x_{in_2}$ are the solutions of the quartic equation
 \begin{align}
 \left.
	 \begin{array}{l l}
	b^2(1+m^2)^2x^4+4k_{3_c}x^3+4k_{2_c}x^2+4k_{1_c}x+ k_{0_c} =0
	\end{array}  \right.   \label{eq:COpair}
\end{align}
which satisfy 
 \begin{align}
 \left.
	 \begin{array}{l l}
	\{ x_{in_1} \ | \ O(x_{in_i}) =0  \}
	\end{array}  \right.   \label{eq:Innerpair}
\end{align}
for $i=1,2$.
 The coefficients of \eqref{eq:COpair} are given by
  \begin{align}
 \left.
	 \begin{array}{l l}
	k_{3_c} = b(1+m^2)(bmn-x'_1) \\
	k_{2_c} = \frac{b^2n^2(3m^2+1)}{2} -b(2mnx'_1+\frac{\eta(1+m^2)}{2}) \\
	\qquad \ \ +x_1^{'2} -\beta^2\rho^2(1+m^2)  \\
	k_{1_c} = b^2mn^3 - bn(nx'_1 + m\eta) +x'_1\eta -2\beta^2\rho^2mn \\
	k_{0_c} = (bn^2-\eta)^2 - (2\beta\rho n)^2 
	\end{array}  \right.   \label{eq:COpairCoef}
\end{align}
 where $b=(1-\beta^2)$ and $\eta=\rho^2-x_1^{'2}$.
 
\end{theorem}

\textit{Proof}. The proof is provided in Appendix A.

 It is important to note that in the presence of two interceptors, $R_1$ and $R_2$, simultaneous blocking of $B_L$ by both interceptors is not always the optimal strategy. If \eqref{eq:CO1} and  \eqref{eq:CO2} do not intersect, then simultaneous blocking is clearly not the optimal strategy. 
Furthermore, intersection of \eqref{eq:CO1} and  \eqref{eq:CO2} is a necessary but not a sufficient condition for simultaneous blocking of $B_L$ by $R_1$ and $R_2$. 

%

In order to determine whether simultaneous capture or individual capture by one of the red interceptors is the optimal strategy, we first need to determine the optimal strategies in the one-on-one case. Without loss of generality consider $R_1$ to be the participating red interceptor in the one-on-one case. The obtained solution can be applied, separately, to each one of the interceptors in order to determine the individual optimal strategies and make the appropriate comparisons to determine the overall optimal strategy. We consider the same relative coordinate frame as in the previous theorem.
Let $\zeta=sign(y'_s)$. Also, let  $c_\phi=\cos\phi$ and $s_\phi=\sin\phi$ where $\phi$ is the Line-of-sight (LOS) angle from $B_L$ to $R_s$ in the relative frame, that is, $\phi=\lambda_s-\lambda_1$.

\begin{theorem}  \label{th:solo}
Consider the attack stage of the BVR differential game and assume that $\textbf{x}\in \mathcal{R}_R$. 
 The optimal headings of players $B_L$ and $R_1$ in the relative coordinate frame are given by the state-feedback policies
 \begin{align}
 \left.
	 \begin{array}{l l}
	 \cos\theta_B^* &= \frac{x^*}{\sqrt{x^{*2} + y^{*2}}}      \\
	    \sin\theta_B^* &= \frac{y^*}{\sqrt{x^{*2} + y^{*2}}}      \\
	     \cos\theta_1^* &= \frac{x^*-x'_1}{\sqrt{(x^*-x'_1)^2 + y^{*2}}}      \\
	          \sin\theta_1^* &= \frac{y^*}{\sqrt{(x^*-x'_1)^2 + y^{*2}}}   
\end{array}  \right.    \label{eq:OptimalInputSolo}
\end{align}
 where $x^*= r^*\cos\theta^*$, $y^*= \zeta r^*\sin\theta^*$,  $\cos\theta^* = \frac{br^{*2} - 2\beta\rho r^* -\eta}{2x'_1r^*}$,  and  $\sin\theta^* = \frac{\sqrt{4x_1^{'2}r^{*2} - [br^{*2} - 2\beta\rho r^* -\eta]^2}}{2x'_1r^*}$.
 The optimal radius, $r^*$, is the solution of the sixth-order equation
 \begin{align}
 \left.
	 \begin{array}{l l}
	k_6 r^6 +k_5r^5+k_4r^4+k_3r^3+k_2r^2+k_1r+ (\beta\rho\eta)^2 =0
	\end{array}  \right.   \label{eq:COsolo}
\end{align}
that minimizes the cost
 \begin{align}
 \left.
	 \begin{array}{l l}
	J\!\!\!&= d_s^2+r^2- \frac{d_s}{x'_1}\big(c_\phi[(1\!-\!\beta^2)r^{*2}\! -\! 2\beta\rho r^* \!-\eta] \\
         &~~ +s_\phi \sqrt{4x_1^{'2}r^2\!-\![(1\!-\!\beta^2)r^{*2}\! -\! 2\beta\rho r^* \!-\eta]^2}\big)
	\end{array}  \right.    \label{eq:costSolo}
\end{align}
where
 \begin{align}
 \left.
	 \begin{array}{l l}
	k_1 = 4\beta\rho[\eta(\beta^2\rho^2\!-\!x_1^{'2}s_\phi^2 )+\frac{1}{2}\eta^2(\frac{x'_1}{d_s}c_\phi\!-\!b)]  \\
	k_2 = 2\beta^2\rho^2[2 \beta^2\rho^2\!-\! 2x_1^{'2}(1\!+\!s_\phi^2) \!+\! \eta(4\frac{x'_1}{d_s}c_\phi\!-\!5b) ] \\
	\qquad +\eta^2[b^2+\frac{x'_1}{d_s}(\frac{x'_1}{d_s}\!-\!2p)] \! +\! 4x_1^{'2}s_\phi^2(x_1^{'2}\!+\!b\eta) \\
	k_3 = 4\beta\rho\big( \beta^2\rho^2(2\frac{x'_1}{d_s}c_\phi\!-\!3b) \!+\! b\eta(2b-3\frac{x'_1}{d_s}c_\phi) \\
	\qquad \ \ +  x_1^{'2}[ \frac{\eta}{d_s^2} +b(2\!+\!s_\phi^2) -2\frac{x'_1}{d_s}c_\phi ]\big)  \\
	k_4 = \beta^2\rho^2[13b^2\!+\!4\frac{x'_1}{d_s}(\frac{x'_1}{d_s}\!-\!4p)] \\
	\qquad -2(b\eta\!+\!2x_1^{'2})[b^2+\frac{x'_1}{d_s}(\frac{x'_1}{d_s}\!-\!2p)]  \\
	k_5 = -2b\beta\rho[3b^2+\frac{x'_1}{d_s}(2\frac{x'_1}{d_s}-5p)]  \\
	k_6 = b^2[b^2+\frac{x'_1}{d_s}(\frac{x'_1}{d_s}-2p)].
	\end{array}  \right.   \label{eq:COcoeffSolo} 
\end{align}
Also, $b=(1-\beta^2)$, $\eta=\rho^2-x_1^{'2}$, and  $p=bc_\phi$.
\end{theorem}
\textit{Proof}.  The proof is provided in Appendix A.

\textit{Determining the optimal strategy}. Based on the results of the previous two theorems we now show how to determine the optimal strategy in the case where \eqref{eq:CO1} and  \eqref{eq:CO2} intersect each other. The dominance region of $\mathcal{B}_L$ is first obtained.  Let $(I_x^*,I_y^*)$ denote the cooperative intersection point; this intersection point is obtained from Theorem \ref{th:coop}.
Then, the individual solutions are computed. Let $(x_1^*,y_1^*)$ and $(x_2^*,y_2^*)$ denote the individual aimpoints with respect to $R_1$ and $R_2$ (an additional change of coordinates is necessary to obtain $(x_2^*,y_2^*)$); these aimpoints are obtained from Theorem \ref{th:solo}.
Three cases may occur. Case 1) If $(x_i^*,y_i^*)\notin \mathcal{B}_L$ for both $i=1,2$, then, $(I_x^*,I_y^*)$ is the optimal solution. Case 2) If $(x_i^*,y_i^*)\in \mathcal{B}_L$ for only one $i=1$ or $i=2$, then, $(x_i^*,y_i^*)$ is the optimal strategy. 
Case 3) If $(x_i^*,y_i^*)\in \mathcal{B}_L$ for both $i=1$ and $i=2$, then, the optimal strategy is given by $(x^*,y^*)$ such that $y^*=mx^*+n$ and $x^*=\arg\min \{J(x^*_1), J(x^*_2) \}$.
\section{The Retreat Stage} \label{sec:Retreat}

\subsection{Problem Formulation} \label{sec:RetreatFor}

In this section we consider the retreat stage which occurs after $B_L$ is engaged by the interceptors.  
 In the retreat stage, $B_L$ is being pursued by two attacking missiles $A_1$ and $A_2$. $B_W$ assists $B_L$ by firing two defending missiles $D_1$ and $D_2$ in order to intercept the attacking missiles and protect the asset $B_L$. The differential game of protecting a valuable asset has been extensively addressed in the case of a cooperating asset \cite{Garcia2019}. However, we are now faced with not one but two pairs of attacking and defending missiles and it is unclear how the asset is able to cooperate with both $D_1$ and $D_2$. In order to solve this problem we consider first the differential game of a non-cooperative asset; in such a case, $B_L$ executes a constant flight path route while a pair of attacking and defending missiles play a differential game of attacking/defending the asset $B_L$. The same differential game can be played by more pairs of attacking and defending missiles. Then, a cost metric is devised in order to consider all missiles and derive an optimal heading for $B_L$ which employs the solution of non-cooperative asset game as building block.
 
 Let us consider the non-cooperative asset $B_L=(x_B,y_B)$, an attacking missile $A=(x_A,y_A)$ and a defending missile $D=(x_D,y_D)$ in the fixed frame. $B_L's$ constant heading is denoted by $\theta_B$ and it is known to both $A$ and $D$.
 The speeds of the players are constant they are denoted by $v_B$, $v_A$, and $v_D$, respectively. The agents have simple motion as it is commonly found in the games of Isaacs \cite{Isaacs65}.
The complete state of the game is defined by $\textbf{x}:=( x_B, y_B, x_A, y_A, x_D, y_D)\in \mathbb{R}^6$. 
We note that, although $B_L$'s heading is constant during the retreat stage, its state is still needed by the remaining players to determine their optimal strategies.
The Attacker's control variable is his instantaneous heading angle, $\textbf{u}_A=\left\{\chi\right\}$. The Defender's control variable is $\textbf{u}_D=\left\{\psi\right\}$. The dynamics $\dot{\textbf{x}}=\textbf{f}(\textbf{x},\textbf{u}_A,\textbf{u}_D)$ are defined by the system of ordinary differential equations
\begin{align}
 \left.
	\begin{array}{l l}
        \dot{x}_B&=v_B\cos\theta_B,    \ \ \ \ \ x_B(0)=x_{B_0}   \\
	\dot{y}_B&=v_B\sin\theta_B,  \ \ \ \ \ y_B(0)=y_{B_0}  \\
	\dot{x}_A&=v_A\cos\chi,   \ \ \ \ \ \ x_A(0)=x_{A_0}  \\
	\dot{y}_A&=v_A\sin\chi,   \ \ \ \ \ \ y_A(0)=y_{A_0}  \\
	\dot{x}_D&=v_D\cos\psi,   \ \ \ \ \ \ x_D(0)=x_{D_0}  \\
	\dot{y}_D&=v_D\sin\psi,   \ \ \ \ \ \ y_D(0)=y_{D_0}
	\end{array}  \right.   \label{eq:DynRetr}
\end{align}
where the admissible controls are given by $\chi,\psi \in [-\pi,\pi)$. 
 We assume that the missiles $A$ and $D$ have the same speed while $B_L$ is slower than the missiles. Define the speed ratio $\alpha=v_B/v_A<1$. 
The initial state of the system is defined as
\begin{align}
	\textbf{x}_0 := (x_{B_0}, y_{B_0}, x_{A_0}, y_{A_0}, x_{D_0}, y_{D_0}) = \textbf{x}(t_0).  \nonumber
\end{align}
The termination set is defined as follows
\begin{align}
   \mathcal{T} :=  \mathcal{T}_A    \   \bigcup \   \mathcal{T}_D   \label{eq:TwoSetsRetr}
\end{align}
where 
\begin{align}
   \mathcal{T}_A=   \big\{ \ \textbf{x} \ |  \sqrt{ (x_B-x_A)^2+(y_B-y_A)^2} = 0 \big\}   \label{eq:Set1Retr}
\end{align}
represents the outcome where $B_L$ is captured by $A$. On the other hand
\begin{align}
   \left.
	\begin{array}{l l}
	 \mathcal{T}_D:= \big\{ \ \textbf{x} \ | \sqrt{ (x_D-x_A)^2 + (y_D-y_A)^2} = 0 \big\}   \label{eq:Set2Retr}
\end{array}  \right. 
\end{align}
represents the outcome where $A$ is intercepted by $D$ before $A$ can capture $B_L$. 

In this section we focus on the case where $A$ is successfully intercepted by $D$. The opposing case will be addressed in future research. Let $\mathcal{R}_D\in \mathbb{R}^6$ denote the Defender's winning subspace where, under optimal play, $A$ is successfully intercepted by $D$.
The terminal time $t_f$ is defined as the time instant when the state of the system satisfies \eqref{eq:Set2Retr}, at which time the terminal state is $\textbf{x}_f: = (x_{B_f}, y_{B_f}, x_{A_f}, y_{A_f}, x_{D_f}, y_{D_f}) = \textbf{x}(t_f)$.
The terminal cost/payoff functional is
\begin{align}
  J(\textbf{u}_A(t),\textbf{u}_D(t),\textbf{x}_0)=\Phi(\textbf{x}_f)	\label{eq:costDG2Retr}
\end{align}
where
\begin{align}
  \Phi(\textbf{x}_f):=\sqrt{(x_{A_f}-x_{B_f})^2+(y_{A_f}-y_{B_f})^2}.	\label{eq:costDGRetr}
\end{align}
The cost/payoff functional depends only on the terminal state - the ATDDG is a terminal cost/Mayer type game. Its Value is given by
\begin{align}
  V(\textbf{x}_0):= \min_{\textbf{u}_A(\cdot)} \ \max_{\textbf{u}_D(\cdot)} J(\textbf{u}_A(\cdot),\textbf{u}_D(\cdot);\textbf{x}_0)	\label{eq:costDG3Retr}
\end{align}
subject to \eqref{eq:DynRetr} and \eqref{eq:Set2Retr}, where $\textbf{u}_A(\cdot)$ and $\textbf{u}_D(\cdot)$ are the players' state feedback strategies.

 \begin{corollary}  \label{cor:SLretr}
Consider the retreat stage of the BVR differential game and assume that $\textbf{x}\in\mathcal{R}_D$. The headings of the players $A$ and $D$ are constant under optimal play and the optimal trajectories are straight lines.
\end{corollary}

\begin{figure}
	\begin{center}
		\includegraphics[width=8.4cm,trim=1.5cm .9cm 1.4cm .7cm]{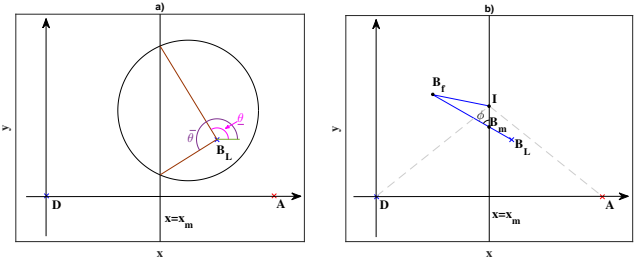}
	\caption{Active defense of non-cooperative asset. a) Game of Kind. b) Derivation of optimal strategies }
	\label{fig:TADfig}
	\end{center}
\end{figure}

\subsection{Optimal Strategies in the Retreat Stage} \label{sec:RetreatStr}

  Without loss of generality we analyze this problem using the relative coordinate frame shown in Fig. \ref{fig:TADfig}. In this coordinate frame we have that $D=(0,0)$, $A=(x'_A,0)$, and $B_L=(x'_B,y'_B)$. The locations of the agents in the relative coordinate frame can be readily obtained from the positions of the agents in the fixed frame as follows
 \begin{align}
\left.
	\begin{array}{l l}
  	x'_A=\sqrt{(x_A-x_D)^2+(y_A-y_D)^2} \\
	x'_B=d_B\cos(\lambda_B-\lambda_A) \\
	y'_B=d_B\sin(\lambda_B-\lambda_A) 
         \end{array}   \right.  \nonumber 
\end{align}
where $d_B=\sqrt{(x_B-x_D)^2+(y_B-y_D)^2}$, $\lambda_B=\arctan(\frac{y_B-y_D}{x_B-x_D})$, and $\lambda_A=\arctan(\frac{y_A-y_D}{x_A-x_D})$.  Also define $x_m=\frac{1}{2}x'_A$. 
 
\begin{theorem}  \label{th:GoKreteat}
Consider the retreat stage of the BVR differential game and assume that $\textbf{x}\in\mathcal{R}_D$. The Defender is able to intercept the Attacker before the latter captures $B_L$ if $\theta_B$ is such that
\begin{align}
\left.
	\begin{array}{l l}
  	\theta_B  \in [\underline{\theta}+\lambda_A,\bar{\theta}+\lambda_A]
         \end{array}   \right.   \label{eq:RelFgokTH}
\end{align}
where $\underline{\theta}=\arctan(\frac{\underline{y}-y'_B}{x_m-x'_B})$, $\bar{\theta}=\arctan(\frac{\bar{y}-y'_B}{x_m-x'_B})$, and
\begin{align}
\left.
	\begin{array}{l l}
  	\underline{y} = \frac{y'_B+\sqrt{\alpha^2[(x'_B-x'_A)^2+y_B^{'2}] - [(1\!-\!\alpha^2)x_m-x'_B+\alpha^2 x'_A]^2}}{1-\alpha^2} \\
        \bar{y} = \frac{y'_B-\sqrt{\alpha^2[(x'_B-x'_A)^2+y_B^{'2}] - [(1\!-\!\alpha^2)x_m-x'_B+\alpha^2 x'_A]^2}}{1-\alpha^2}
        \end{array}   \right.  \nonumber 
\end{align}
\end{theorem}
\textit{Proof}.  The proof is provided in Appendix B.

Since we assume that $\textbf{x}\in\mathcal{R}_D$, there exist a non-empty set of heading angles such that $B_L$ is able to escape. In other words, the Apollonius circle between $A$ and $B_L$ intersects the line $x=x_m$. Then, the Game of Degree is played in the winning region of the asset/defender team for some fixed $\theta_B$ such that \eqref{eq:RelFgokTH} holds, where the objective functional is \eqref{eq:costDG2Retr}.

\begin{theorem}  \label{th:TADnonC}
Consider the retreat stage of the BVR differential game. Assume that $B_L$ moves with constant heading $\theta_B$. The state-feedback optimal strategies of $A$ and $D$ in the relative coordinate frame are given, respectively, by
  \begin{align}
\left.
	\begin{array}{l l}
  	 \cos\chi^* &= \frac{x_m-x'_A}{\sqrt{(x_m-x'_A)^2 + y^{*2}}}      \\
	    \sin\chi^* &= \frac{y^*}{\sqrt{(x_m-x'_A)^2 + y^{*2}}}      \\
	     \cos\psi^* &= \frac{x_m}{\sqrt{x_m^2 + y^{*2}}}      \\
	          \sin\psi^* &= \frac{y^*}{\sqrt{x_m^2 + y^{*2}}} 
         \end{array}   \right.  \label{eq:NonCoopOpt}
\end{align}
where  $y^*$ is the solution of the following quartic equation
\begin{align}
\left.
	\begin{array}{l l}
  	c_4y^4 +2c_3y^3+c_2y^2+2c_1y+c_0=0 
	 \end{array}   \right.  \label{eq:NonCoopQeq}
\end{align}
which minimizes the cost 
\begin{align}
\left.
	\begin{array}{l l}
  	J=(d_m+\alpha\sqrt{x_m^2+y^2})^2 +(y-y_m)^2 \\
	\qquad -2(y-y_m)(d_m+\alpha\sqrt{x_m^2+y^2})\cos\varphi.
	 \end{array}   \right.  \label{eq:NonCoopCost}
\end{align}
The coefficients of \eqref{eq:NonCoopQeq} are given by
\begin{align}
\left.
	\begin{array}{l l}
  	c_4 =(2\alpha\cos\varphi)^2-(1+\alpha^2)^2 \\
	c_3= y'_B+\alpha^2[ (1\!-\!2\cos^2\varphi)y_m -d_m\cos\varphi ] \\
	c_2= \alpha^2[(d_m\!+\!y_m\cos\varphi)^2 +(2x_m\cos\varphi)^2] \\
	\qquad-y_B^{'2} - x_m^2(1+\alpha^2)^2\\
	c_1= y'_Bx_m^2 + \alpha^2x_m^2y_m\sin^2\phi) \\
	c_0= (\alpha x_m^2\cos\varphi)^2-x_m^2y_B^{'2}
	\end{array}   \right.  \nonumber 
\end{align}
where $\varphi=\theta_B-\frac{\pi}{2}-\lambda_A$, $d_m=\frac{x_m-x'_B}{\sin\varphi}$, and $y_m=y'_B-d_m\cos\varphi$.
\end{theorem} 
\textit{Proof}.  The proof is provided in Appendix B.

\subsection{Cooperative Strategy of Blue Leader} \label{sec:RetreatBLstr}
Given a heading that guarantees escape, \textit{i.e.} $\theta_B$ satisfies \eqref{eq:RelFgokTH}, the results in Theorem  \ref{th:GoKreteat} and Theorem  \ref{th:TADnonC} can be used to obtain the optimal strategies of each separate pair of attacking and defending missiles. In addition, $B_L$ can search for an optimal heading that balances off a linear combination of the terminal costs associated to each separate game. This is achieved as explained in the remaining of this section.

Consider $B_L$ being pursued by two attacking missiles $A_1$ and $A_2$. $B_W$ assists $B_L$ by firing two defending missiles $D_1$ and $D_2$. We first apply Theorem  \ref{th:GoKreteat} to each pair of attacking-defending missiles. From the $A_1-D_1$ pair we obtain $\underline{\theta}_1$ and $\bar{\theta}_1$ and from the $A_2-D_2$ pair we obtain $\underline{\theta}_2$ and $\bar{\theta}_2$. Thus, the feasible heading of $B_L$, in order to successfully escape from both $A_1$ and $A_2$, is such that 
\begin{align}
\left.
	\begin{array}{l l}
  	\theta_B  \in [\theta_l,\theta_u]
	 \end{array}   \right.  \label{eq:FeasbThetaB}
\end{align}
 where $\theta_l=\max \{ \underline{\theta}_1+\lambda_{A_1},\underline{\theta}_2+\lambda_{A_2}\}$, $\theta_u=\min\{\bar{\theta}_1+\lambda_{A_1},\bar{\theta}_2+\lambda_{A_2}\}$, $\lambda_{A_1}=\arctan(\frac{y_{A_1}-y_{D_1}}{x_{A_1}-x_{D_1}})$, and $\lambda_{A_2}=\arctan(\frac{y_{A_2}-y_{D_2}}{x_{A_2}-x_{D_2}})$. 

In order to optimally evade both $A_1$ and $A_2$, we define the objective functional  
\begin{align}
\left.
	\begin{array}{l l}
  	J_c (\theta_B) =  w V_1(\textbf{x}_1,\theta_B) 
	 + (1-w) V_2(\textbf{x}_2,\theta_B)
	 \end{array}   \right.  \label{eq:CompCost}
\end{align}
for parameter  $w\in (0,1)$ and subject to $\theta_B  \in [\theta_l,\theta_u]$. Also, $\textbf{x}_1=( x_B, y_B, x_{A_1}, y_{A_1}, x_{D_1}, y_{D_1})$ and $\textbf{x}_2=( x_B, y_B, x_{A_2}, y_{A_2}, x_{D_2}, y_{D_2})$. Clearly, $B_L$ tries to maximize \eqref{eq:CompCost}.

In order to obtain $V_i(\textbf{x}_i,\theta_B)$, for $i=1,2$, we apply the results of  Theorem  \ref{th:TADnonC}, that is, $V_i(\textbf{x}_i,\theta_B)$ is the Value of the game played by $A_i$ and $D_i$ given some $\theta_B\in [\theta_l,\theta_u]$.  Hence, a numerical search can be performed over the admissible set of heading values $\theta_B\in [\theta_l,\theta_u]$ in order to maximize \eqref{eq:CompCost}. In summary, the results of the differential games, one for each pair of attacking-defending missiles, are used by $B_L$ in order to cooperate with both defenders and maximize a linear combination of both terminal distances.

\section{Examples} \label{sec:Ex}
Consider the initial positions $B_{L_0}=(-6, 8)$,  $B_{W_0}=(-7.4,6.2)$,   $R_{1_0} =(15, 14)$, and  $R_{2_0}=(16, 6.5)$. The red stationary asset is located at $R_s=(15.5, 10)$. The air-to-air firing range is $\rho=5.0$ and the air-to-ground firing range is $\rho_s=7$. During the attack stage we have that $v_B=1$ and $v_1=v_2=1.25$; hence, the speed ratio parameter is $\beta=1.25$. Solving the attack stage using simple motion dynamics we obtain that the Value of the game in this stage is $V=13.9870$. 

\begin{figure}
	\begin{center}
		\includegraphics[width=8.0cm,trim=1.1cm .9cm 1.2cm .7cm]{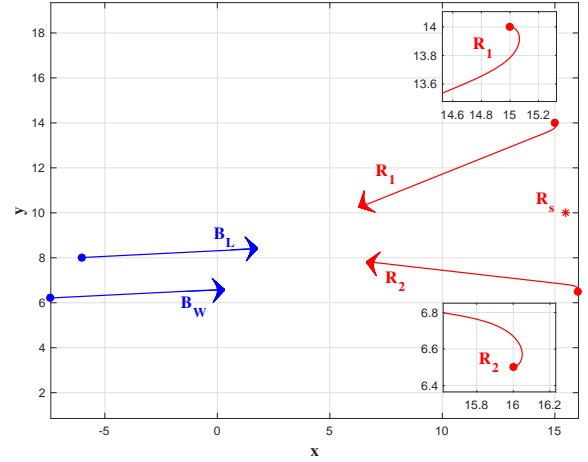}
	\caption{Attack Stage. Trajectories with zoom in view about the initial paths of the red interceptors}
	\label{fig:Ex1attack}
	\end{center}
\end{figure}

\begin{figure}
	\begin{center}
		\includegraphics[width=8.0cm,trim=1.1cm .9cm 1.2cm .7cm]{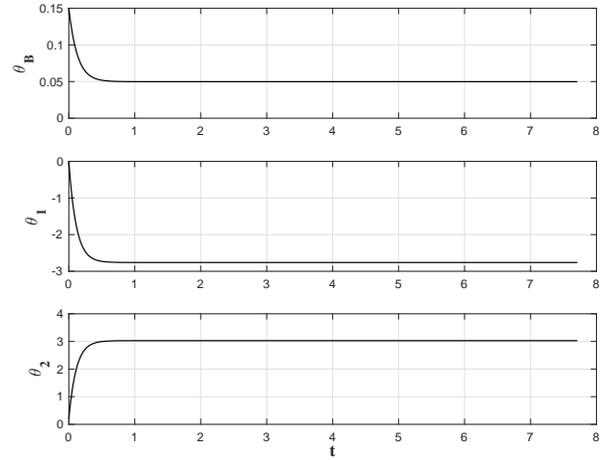}
	\caption{Attack Stage. Heading profiles}
	\label{fig:Ex1headings}
	\end{center}
\end{figure}

We are interested in applying these results to a more realistic case where the vehicles exhibit turning rate constraints. We employ the typical first-order model dynamics \cite{Pachter94} which is described as follows.
\begin{align}
 \left.
	\begin{array}{l l}
      \dot{x}_i&=v_i\cos\theta_i,  \qquad \qquad x_i(0)=x_{i_0} \\
        \dot{y}_i&= v_i \sin\theta_i, \qquad \qquad  y_i(0)=y_{i_0}   \\
          \dot{\theta}_i&=\frac{1}{\tau_i}(\theta_{c_i}-\theta_i), \qquad \  \theta_i(0)=\theta_{i_0}.  
	\end{array}  \right.   \label{eq:DynamicsL}
\end{align}
for $i=B,1,2$. We consider the parameters $\tau_B=0.14$ and $\tau_{1}=\tau_{2}=0.12$.
The procedure to apply the results in this paper to the models in \eqref{eq:DynamicsL} is as follows. The optimal strategies obtained in this paper return the optimal headings $\theta_i^*$ for each player $i=B,1,2$. We now make $\theta_{c_i}(t)=\theta_i^*(t)$ in order to guide each player to acquire its own optimal heading. The optimal heading is expected to be time-varying, even when all players implement their optimal strategies, since the players are now constrained on how fast they can turn. At each time instant, each player uses the current state of the game in order to update its optimal control heading.
Fig. \ref{fig:Ex1attack} shows the trajectories during the attack stage where players $B_L$, $R_1$, and $R_2$ aim at the optimal blocking point which is given by the intersection of the two Cartesian Ovals. However, at the start of the engagement, the red interceptors have a significant initial heading disadvantage with respect to the blue UAVs. This is illustrated in the zoom in plots within Fig. \ref{fig:Ex1attack}. In addition, the heading profiles of $B_L$, $R_1$ and $R_2$ are shown in Fig. \ref{fig:Ex1headings} where it can be seen that $R_1$ and $R_2$ take a harder turn than $B_L$ at the start of the engagement. The red interceptors are nevertheless able to simultaneously block $B_L$ and the terminal separation between $B_L$ and $R_s$ is $J=13.9052$ which is just slightly smaller than the Value of the game. The effects of turning rate constraint are negligible at the BVR range.


\begin{figure}
	\begin{center}
		\includegraphics[width=8.0cm,trim=1.1cm .9cm 1.2cm .2cm]{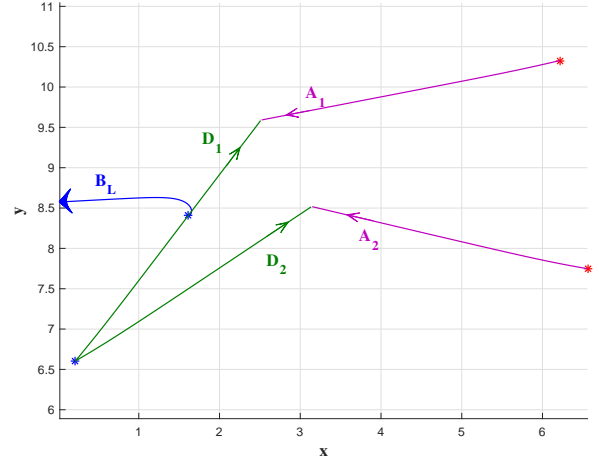}
	\caption{Retreat Stage. Trajectories where $B_L$ is subject to turning rate and acceleration constraints}
	\label{fig:Ex1retr}
	\end{center}
\end{figure}

We now consider the retreat stage of the differential game. At the end of the attack stage, both interceptors reach a distance $\rho$ within $B_L$ simultaneously and  $B_L$  launches an attack against the interceptors and retreats. 
 $R_1$ and $R_2$ fire one missile each against $B_L$ from their position at the end of the attack stage; these missiles are denoted by $A_1$ and $A_2$ in Fig. \ref{fig:Ex1retr}. $B_W$ fires two missiles to protect $B_L$ from its position at the end of the attack stage which are denoted by $D_1$ and $D_2$. The positions of the players at the end of the attack stage which is also the beginning of the retreat stage are denoted by $*$ marks. Note that the retreat stage takes place at a significant closer range than the attack stage.
 
 For the retreat stage, $B_L$ includes the sense of urgency and retreats at higher speed. $B_L$ increases its speed from $v_B=1$ in the attack stage to $v_B=1.5$ in the retreat stage. The missiles speed in this example is $v_A=v_D=3$ and the speed ratio parameter, when $B_L$ flies at max speed, is $\alpha=0.5$.  In order to avoid the attacking missiles flying near its teammate, $B_W$, we include an additional constraint on the optimal heading of $B_L$ during the retreat stage. This is done by restricting $B_L$ to turn in the direction opposite to the location of $B_W$. Then, in this example, the additional constraint is $\theta_B<\theta_{B_f}^a+\pi$, where $\theta_{B_f}^a$ is the heading of $B_L$ at the end of the attack stage. In this example, the optimal heading of $B_L$ is $\theta^*_B=3.1956$ which is equal to the imposed constraint. The Value of the game in the retreat stage, for parameter $w=0.5$ and assuming all vehicles have simple motion, is $V=3.2442$. 
 
\begin{figure}
	\begin{center}
		\includegraphics[width=7.9cm,trim=1.1cm .9cm 1.2cm .0cm]{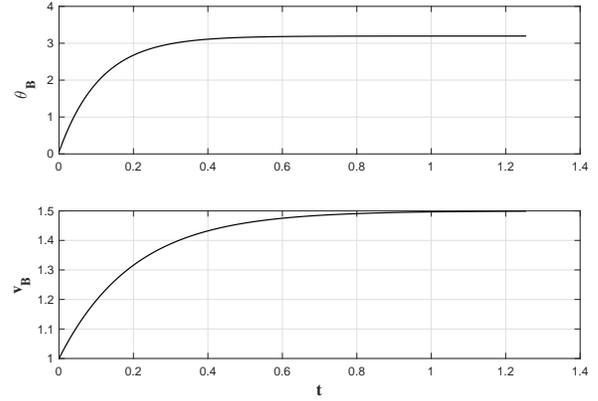}
	\caption{Retreat Stage. Heading and speed profiles of $B_L$}
	\label{fig:Ex1retrProf}
	\end{center}
\end{figure}

 In the trajectories shown in Fig. \ref{fig:Ex1retr} we apply the results obtained in Section \ref{sec:Retreat} to the case where $B_L$ has both turning rate constraints and acceleration constraints. The attacker missiles are intercepted at different times, $A_2$ is intercepted about 0.1 seconds earlier than $A_1$. Fig. \ref{fig:Ex1retr} shows the complete trajectory of $B_L$, until $A_1$ is intercepted. $B_L$'s speed profile is given by the first-order dynamical model $\dot{v}_B=\frac{1}{\tau_{B_v}}(v_{c_B}-v_B)$ where $\tau_{B_v}=0.2$. The cost in this case where turning rate and acceleration constraints are considered is $J_c=2.9308$ which is smaller than the Value of the game. The heading and speed profile of $B_L$ are shown in Fig. \ref{fig:Ex1retrProf}. 
 
 The duration of the retreat stage is much shorter than that of the attack stage and the engagement occurs at closer range than the engagement at the attack stage. Additionally, $B_L$ is not only subject to turning rate constraints but also to acceleration constraints. These are the main reasons for the decrease in performance by $B_L$. The approach and the strategies obtained assuming simple motion models work well for the retreat stage and they are robust with respect to constraints which are not modeled in the problem formulation. The defenders are able to react to the variations introduced by the dynamic constraints of $B_L$ and intercept the attackers  by implementing the obtained state-feedback strategies. The performance variation, although smaller than $10\%$, is not negligible as it was in the attack stage. 
 
 In future work we plan to address the retreat stage, with $B_L$'s turning rate and acceleration constraints, in more detail. We note that a model of a player with four states, with heading and speed included as states, significantly complicates the analysis. Only open-loop numerical solutions are usually obtained when considering this kind of model. Since open-loop solutions are not useful for implementation in the type of attack and defense problem under consideration. We plan to approach the problem by analyzing the heading and speed profiles of the players and updating the state-feedback strategies according to the current values of the parameters. For instance, the current speed of $B_L$ can be used to update the speed ratio parameter and, in turn, update the state-feedback strategies of the players. Since the Value of the game obtained while assuming simple motion models is not achievable, $B_L$ cannot instantaneously acquire its max speed, the proposed analysis will provide a better approximation of the achievable value of the game and a refinement of the optimal strategies.



\section{Conclusion} \label{sec:concl}
 This paper formulated an operational relevant combat scenario in the Beyond Visual Range as a two-team, zero-sum differential game. This game consists of two stages: the attack stage and the retreat stage. The optimal strategies of each stage were obtained in this paper. In the attack stage, these strategies leveraged cooperative actions between agents of the same team in order to prevent an opposing player from reaching the engagement zone of the protected asset and minimizing risk by blocking enemies the farthest away from the asset. In the retreat stage, the active target defense paradigm was used and extended in order consider two pairs of attacking and defending missiles and for the target to cooperate with its two defenders.

\section*{Appendix A. Proofs of Attack Stage} \label{sec:App}  
\textbf{Proof of Theorem \ref{th:SL}}. 
The Hamiltonian of the differential game is
\begin{align}
  \left.
	\begin{array}{l l}
	\mathcal{H}\!\!\!&=\varkappa_{x_B}\cos\theta_B + \varkappa_{y_B}\sin\theta_B +\beta\varkappa_{x_1}\cos\theta_1 \\ 
	&~~ + \beta\varkappa_{y_1}\sin\theta_1 +\beta\varkappa_{x_2}\cos\theta_2 + \beta\varkappa_{y_2}\sin\theta_2 
\end{array}  \right. 
\end{align}
where $\varkappa=(\varkappa_{x_B},\varkappa_{y_B},\varkappa_{x_1},\varkappa_{y_1},\varkappa_{x_2},\varkappa_{y_2})$ is the vector of co-states. Note that the Hamiltonian and the dynamics are separable (or decoupled) in the controls $\theta_B$ and $\theta_R$. Hence, $\min_{\theta_B} \max_{\theta_1,\theta_2}  \mathcal{H}= \max_{\theta_1,\theta_2} \min_{\theta_B}  \mathcal{H}$ and Isaacs' condition holds.

The optimal control inputs in terms of the co-state variables are obtained from Isaacs' Main Equation 1 (ME 1)
\begin{align}
   \min_{\theta_B} \max_{\theta_1,\theta_2} \mathcal{H} =0 	\label{eq:minmax}
\end{align}
and they are characterized by the relationships
\begin{align}
\left.
	\begin{array}{l l}
  &\cos\theta_B^*=-\frac{\varkappa_{x_B}}{\sqrt{\varkappa_{x_B}^2+\varkappa_{y_B}^2}}, \ \ \ \  \sin\theta_B^*=-\frac{\varkappa_{y_B}}{\sqrt{\varkappa_{x_B}^2+\varkappa_{y_B}^2}}  \\
	& \cos\theta_1^*=\frac{\varkappa_{x_1}}{\sqrt{\varkappa_{x_1}^2+\varkappa_{y_1}^2}}, \ \ \ \ \ \ \   \sin\theta_1^*=\frac{\varkappa_{y_1}}{\sqrt{\varkappa_{x_1}^2+\varkappa_{y_1}^2}}  \\
	& \cos\theta_2^*=\frac{\varkappa_{x_2}}{\sqrt{\varkappa_{x_2}^2+\varkappa_{y_2}^2}}, \ \ \ \ \ \ \   \sin\theta_2^*=\frac{\varkappa_{y_2}}{\sqrt{\varkappa_{x_2}^2+\varkappa_{y_2}^2}}
	\end{array}  \right.  \label{eq:psico}  
 \end{align}
 The co-state dynamics are obtained from $\dot{\varkappa}=-\frac{\partial H}{\partial \textbf{x}}$ which results in: $\dot{\varkappa}_{x_B}=\dot{\varkappa}_{y_B}=\dot{\varkappa}_{x_1}=\dot{\varkappa}_{y_1}=\dot{\varkappa}_{x_2}=\dot{\varkappa}_{y_2}=0$; hence, all co-states are constant and we have that the optimal headings $\theta_B^*$, $\theta^*_1$ and $\theta_2^*$ are constant as well. \qquad \qquad \quad $\square$

\textbf{Proof of Theorem \ref{th:coop}}. 
In order to determine the intersection points of the COs \eqref{eq:CO1} and \eqref{eq:CO2}, we substitute the generic coordinates $(x,y)$ by $(I_x,I_y)$, the coordinates of the intersection points into these two equations. Hence, we have that 
\begin{align}
 \left.
	 \begin{array}{l l}
	  \sqrt{(I_x-x'_1)^2+I_y^2} =  \sqrt{(I_x-x'_2)^2+(I_y-y'_2)^2}   
\end{array}  \right.   \nonumber
\end{align}
that is, the distance between the current position of $R_1$ and the intersection point $I$ is the same as the distance between the current position of $R_2$ and the intersection point $I$. Therefore, the intersection point lies on the orthogonal bisector of the segment $\overline{R_1R_2}$ which is defined by $y=mx+n$ where $m$ and $n$ are given by \eqref{eq:mni}. 

It is possible to obtain the intersections of $y=mx+n$ with either  \eqref{eq:CO1} or \eqref{eq:CO2}. Without loss of generality, consider  \eqref{eq:CO1}. Taking the square of both sides of \eqref{eq:CO1} and rearranging terms we obtain
\begin{align}
 \left.
	 \begin{array}{l l}
	 b(x^2+y^2)  - 2x'_1x -\eta = 2\beta\rho\sqrt{x^2+y^2}.
\end{array}  \right.   \nonumber
\end{align}
Once again, taking the square of both sides of the previous equation and combining like terms we arrive at
\begin{align}
 \left.
	 \begin{array}{l l}
	 b^2(y^4 + x^4) - 4bx'_1x^3 + 4[\frac{b^2x^2}{2} - b(x'_1x + \frac{\eta}{2}) -\beta^2\rho^2] y^2 \\
	 + 4(x_1^{'2} -\frac{b\eta}{2} -\beta^2\rho^2)x^2 +4x'_1\eta x +\eta^2 =0.
\end{array}  \right.   \nonumber
\end{align}
Now, substitute $y=mx+n$ into the previous equation. Combining like terms, we obtain the quartic equation \eqref{eq:COpair} with coefficients given by \eqref{eq:COpairCoef}.

By taking the square of the previous equations, additional ovals are introduced. In general, each Cartesian Oval consists of two ovals: the inner oval and the outer oval. However, for each case, \eqref{eq:CO1} and \eqref{eq:CO2}, the applicable oval is the inner oval. This is due to the fact that points on the outer ovals do not satisfy the capture conditions \eqref{eq:CO1} or \eqref{eq:CO2}.

If the four solutions of \eqref{eq:COpair} are real, then, the line $y=mx+n$ intersects both the inner ovals the outer ovals.
However, the outer oval is irrelevant and we discard those solutions by imposing the condition \eqref{eq:Innerpair}. Thus, the solutions $x_{in_1}$ and $x_{in_2}$ correspond to the two intersection points of the inner ovals \eqref{eq:CO1} and \eqref{eq:CO2}. The optimal intersection is the closest to the red asset $R_s$ and the optimal headings are given by \eqref{eq:OptimalInputsPCRcoop}. Finally, the optimal headings in the fixed frame, $\theta_{i_f}^*$, are simply given by $\theta_{i_f}^*=\theta_i^*+\lambda_1$, for $i=B,1,2$.  \qquad $\square$

\begin{figure}
	\begin{center}
		\includegraphics[width=7.8cm,trim=1.5cm 1.0cm 1.4cm .3cm]{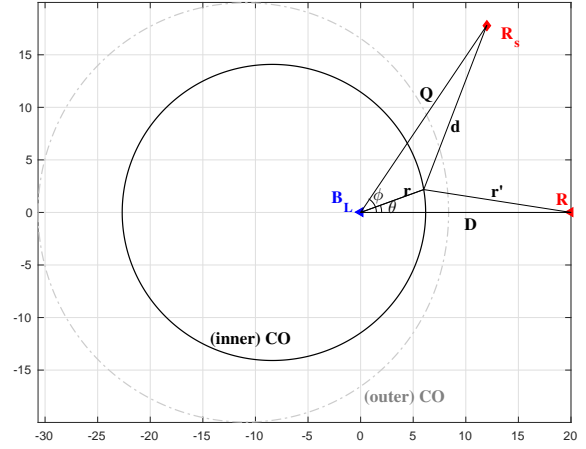}
	\caption{CO between $B_L$ and $R$. Both, inner and outer ovals are shown; however, outer oval does not satisfy the CO condition \eqref{eq:CO1}. $D=x'_1$, $Q=d_s$, $r'=\hat{r}$}
	\label{fig:Indvopt}
	\end{center}
\end{figure}

\textbf{Proof of Theorem \ref{th:solo}}.  
Consider the definition of the Cartesian Oval \eqref{eq:CO1} rewritten in radial form
\begin{align}
\hat{r}=\beta r +\rho.  \label{eq:CObipolar}
\end{align} 
Without loss of generality we consider the configuration in Fig. \ref{fig:Indvopt}.  Using the law of cosines, we have that for any point on the Cartesian Oval
\begin{align}
\left.
	\begin{array}{l l}
        \cos\theta =  \frac{r^2-\hat{r}^2+x_1^{'2}}{2rx'_1}.    
         \end{array}   \right. \nonumber
\end{align}
In addition, $\sin\theta=\frac{1}{2rx'_1}\sqrt{4x_1^{'2}r^2-(r^2-\hat{r}^2+x_1^{'2})^2}$. The distance $d$ is given by
 \begin{align}
\left.
	\begin{array}{l l}
     J(r) = d^2\!\!\! &=d_s^2+r^2-2d_sr\cos(\phi-\theta) \\
         &= d_s^2+r^2- \frac{d_s}{x'_1}\big(\cos\phi(r^2-\hat{r}^2+x_1^{'2}) \\
         &~~ +\sin\phi \sqrt{4x_1^{'2}r^2-(r^2-\hat{r}^2+x_1^{'2})^2}\big).
         \end{array}   \right. \label{eq:SoloCost}
\end{align}
Since $B_L$ is not able to reach $\Gamma$, the optimal strategy is for $B_L$ to be blocked at the closest point with respect to $\Gamma$; since $\Gamma$ is a circular region with center at $R_s$, the closest point to $\Gamma$ is equivalent to the closest point to $R_s$. Such a point can be obtained by taking the derivative of \eqref{eq:SoloCost} with respect to $r$ and setting the result equal to zero as follows
  \begin{align}
\left.
	\begin{array}{l l}
      \frac{d J(r)}{dr}&=2r -\frac{2d_s}{x'_1}\big((r-\beta \hat{r})\cos\phi \\
        &~~+ \sin\phi \frac{2x_1^{'2}r-(r^2-\hat{r}^2+x_1^{'2})(r-\beta \hat{r})}{\sqrt{4x_1^{'2}r^2-(r^2-\hat{r}^2+x_1^{'2})^2}} \big) = 0
         \end{array}   \right. \label{eq:SoloDer}
\end{align}
  where, from \eqref{eq:CObipolar}, $\frac{d\hat{r}}{dr}=\beta$. Equation \eqref{eq:SoloDer} can be written in the following form
   \begin{align}
\left.
	\begin{array}{l l}
   \sin\phi \big(2x_1^{'2}r-(r^2-\hat{r}^2+x_1^{'2})(r-\beta \hat{r}) \big) = \\
     \big((r-\beta \hat{r}) \cos\phi- \frac{x'_1}{d_s}r \big) \sqrt{4x_1^{'2}r^2-(r^2-\hat{r}^2+x_1^{'2})^2}. \\
         \end{array}   \right. \label{eq:SoloDer2}
\end{align}
 Taking the square of both sides of equation  \eqref{eq:SoloDer2} and rearranging terms we obtain the following
  \begin{align}
\left.
	\begin{array}{l l}
    (r^2\!-\!\hat{r}^2\!+\!x_1^{'2})^2 \big((r\!-\!\beta \hat{r})[(1\!-\!\frac{2x'_1}{d_s}\cos\phi)r \!-\!\beta \hat{r}] \!+\! \frac{x_1^{'2}}{d_s^2}r^2  \big) \\
   + 4x_1^{'2}r(r-\beta \hat{r})\big( (\frac{2x'_1}{d_s}\cos\phi\!-\!1)r^2+(\hat{r}^2\!-\!x_1^{'2})\sin^2\phi  \\
    ~~+ r\hat{r} \beta\cos^2\phi  \big) +4x_1^{'4}r^2(\sin^2\phi-\frac{r^2}{d_s^2}) = 0
         \end{array}   \right. \label{eq:SoloDer3}
\end{align}
  In order to simplify the notation, let $c_\phi=\cos\phi$ and $s_\phi=\sin\phi$.
 We substitute \eqref{eq:CObipolar} into \eqref{eq:SoloDer3} in order to write the equation in terms of only one variable, $r$. We have the following
   \begin{align}
\left.
	\begin{array}{l l}
    \Big( (1\!-\!\beta^2)^2r^4 \!-\! 4(1\!-\!\beta^2)\beta\rho r^3 \!+\! 2[2\beta^2\rho^2 \!-\! (1\!-\!\beta^2)\eta]r^2   \\
    ~~ +\! 4\beta\rho\eta r \!+\! \eta^2 \Big) \Big([(1\!-\!\beta^2)^2\!-\! \frac{2x'_1}{d_s}(1\!-\!\beta^2)c_\phi \!+\! \frac{x_1^{'2}}{d_s^2}]r^2  \\
      ~~ + 2\beta \rho[\frac{x'_1}{d_s}c_\phi \!-\! (1\!-\!\beta^2) ]r +\beta^2\rho^2 \Big) \\
      +4x_1^{'2}r\Big( (1\!-\!\beta^2)r\!-\!\beta\rho \Big) \Big( [\frac{2x'_1}{d_s}c_\phi\!-\! (1\!-\!\beta^2) ]r^2  \\
      ~~ + \beta\rho(1\!+\!s_\phi^2)r \!+\! \eta s_\phi^2 \Big)  +4x_1^{'4}r^2(s^2_\phi-\frac{r^2}{d_s^2}) = 0
         \end{array}   \right. \nonumber
\end{align}
 Finally, expanding the terms in the previous equation and combining like terms we arrive at the sixth-order equation \eqref{eq:COsolo} with coefficients given by \eqref{eq:COcoeffSolo}.
  Note that, due to symmetry of the Cartesian Oval about the $x$-axis, two possible points in the Cartesian Oval exist for any given radius $r$: $x=r \cos\theta$, and $y =\pm r\sin\theta$. For the optimal solution, $r^*$, the minimum is achieved on the same side, with respect to the $x$-axis, of the Cartesian Oval as the location of $R_s$; then, we have that  $x^*=r^* \cos\theta^*$, and $y^* =\zeta r^*\sin\theta^*$.

Note that complex roots and real roots of \eqref{eq:COsolo} outside the range of the inner Cartesian Oval $r\in[\underline{r}, \bar{r}]$ can be automatically discarded.  The values $\underline{r}=\frac{x'_1-\rho}{\beta+1}$ and $\bar{r}=\frac{x'_1-\rho}{\beta-1}$ are, respectively, the smallest and largest values of $r$ that satisfy the inner Cartesian Oval equation shown in \eqref{eq:CObipolar}.  

Consider the Cartesian Oval equation, inner and outer oval, in the following form
\begin{align}
\left.
	\begin{array}{l l}
   (1-\beta^2)r^2 - 2(\beta\rho +x'_1\cos\theta)r - \eta=0
         \end{array}   \right. \label{eq:COpolar}
\end{align}
The inner oval is given by
\begin{align}
\left.
	\begin{array}{l l}
   r_{in} (\theta)\!=\! \frac{\beta\rho +x'_1\cos\theta - \sqrt{(\beta\rho +x'_1\cos\theta)^2 + (1-\beta^2)\eta}}{1-\beta^2}
         \end{array}   \right. \label{eq:COinner}
\end{align}
and the outer oval by
\begin{align}
\left.
	\begin{array}{l l}
   r_{out} (\theta)\!= \!\! \frac{\beta\rho +x'_1\cos\theta + \! \sqrt{(\beta\rho +x'_1\cos\theta)^2 +(1-\beta^2)\eta}}{1-\beta^2}
         \end{array}   \right. \label{eq:COouter}
\end{align}
and we obtain $\underline{r}$ when we make $\theta=0$ in \eqref{eq:COinner}. 
 Similarly, the largest value of $r$ is obtained when $\theta=\pi$ in \eqref{eq:COinner}. 
 Thus, $0<\underline{r}\leq r_{in}(\theta) \leq \bar{r}$. 
 
 Furthermore, the roots of \eqref{eq:COsolo} also include the radii corresponding to points in the outer oval (the outer oval is shown in  Fig. \ref{fig:Indvopt} in gray tone) that minimize and maximize the distance with respect to the target point $R_s$. 
We now make $\theta=0$ in \eqref{eq:COouter} and we obtain $r_{out}(\theta=0)= \frac{x'_1+\rho}{1-\beta}<0$. Also, $r_{out}(\theta=\pi)= - \frac{x'_1+\rho}{1+\beta}<0$. In this case, $\frac{x'_1+\rho}{1-\beta}\leq r_{out}(\theta) \leq  -\frac{x'_1+\rho}{1+\beta}<0$. Therefore, $r_{out}<0$ while $r_{in}>0$ and at least two roots of \eqref{eq:COsolo} fall outside the inner oval range, that is, $r_{out}\notin[\underline{r}, \bar{r}]$  and they can be discarded.  \qquad \qquad \quad $\square$

 \section*{Appendix B. Proofs of Retreat Stage}

 \textbf{Proof of Theorem \ref{th:GoKreteat}}. We determine the range of headings of $B_L$ that guarantee successful interception of $A$ by $D$ by studying the dominance regions of the players. Since $A$ is faster than $B_L$, their dominance regions are separated by the Apollonius circle 
 \begin{align}
\left.
	\begin{array}{l l}
  	(x\!-\!\frac{x'_B-\alpha^2 x'_A}{1-\alpha^2})^2 +(y\!-\!\frac{y'_B}{1-\alpha^2})^2 = \alpha^2\frac{(x'_B-x'_A)^2+y_B^{'2}}{1-\alpha^2}
	   \end{array}   \right.  \label{eq:RelFAppCir}
\end{align}
as it is shown in Fig. \ref{fig:TADfig}.a. 
Now, the dominance regions between $A$ and $D$ are separated by the orthogonal bisector of the segment $\overline{AD}$ which is given by the line $x=x_m$. Hence, $B_L$ can escape by reaching the dominance region of $D$; in that way, $D$ is able to aid him by intercepting $A$. The safe segment of the line $x=x_m$ that $B_L$ should aim at is delineated by the intersection of the Apollonius circle \eqref{eq:RelFAppCir} with that line. Substituting $x=x_m$ into \eqref{eq:RelFAppCir} and solving for $y$ we obtain the two intersection points $\underline{y}$ and $\bar{y}$. The corresponding headings in the relative coordinate frame are given by $\underline{\theta}$ and $\bar{\theta}$. 
Finally, accounting for the rotation by $\lambda_A$ used to obtain the relative coordinate frame, we obtain the condition for successful interception of $A$ by $D$ in the fixed frame given in \eqref{eq:RelFgokTH}. \qquad \qquad \qquad \qquad\qquad \qquad\qquad \qquad\qquad \qquad $\square$

 \textbf{Proof of Theorem \ref{th:TADnonC}}. Interception of $A$ by $D$ occurs at a point on the boundary of their dominance regions, which is the line $x=x_m$. Consider the triangle $B_mB_fI$ shown in Fig. \ref{fig:TADfig}.b, where $B_m=(x_m,y_m)$, $B_f$ is the terminal position of $B_L$, and $I=(x_m,y)$ is the interception point in terms of the optimization variable $y$. 
 
 Note that the distance traversed by $B_L$ is equal to $\alpha t_f= \alpha\sqrt{x_m^2+y^2}$. 
The terminal distance between $A$ and $B_L$ is given by
 \begin{align}
\left.
	\begin{array}{l l}
  	J(y)=  \overline{B_mB_f}^2 + (y-y_m)^2 - 2(y-y_m) \overline{B_mB_f}\cos\varphi.
	 \end{array}   \right.  \nonumber 
\end{align}
Additionally, we have that $\overline{B_mB_f}=d_m+\alpha\sqrt{x_m^2+y^2}$ and $J(y)$ can be written as shown in \eqref{eq:NonCoopCost}. Note that $d_m<0$ if $B_L$ is initially located on the $A$-side of the line $x=x_m$ and $d_m>0$ otherwise.

Let us take the derivative of \eqref{eq:NonCoopCost} with respect to $y$ and set it equal to zero as follows
 \begin{align}
\left.
	\begin{array}{l l}
  	\frac{dJ(y)}{dy}= \frac{2\alpha(d_m+\alpha\sqrt{x_m^2+y^2})}{\sqrt{x_m^2+y^2}}y + 2(y-y_m) \\
	 \qquad - 2\cos\varphi[\frac{\alpha(y-y_m)}{\sqrt{x_m^2+y^2}}y +d_m+\alpha\sqrt{x_m^2+y^2}] = 0.
	 \end{array}   \right.  \nonumber 
\end{align}
Multiplying both sides of the previous equation by $\sqrt{x_m^2+y^2}$ and rearranging terms we obtain
  \begin{align}
\left.
	\begin{array}{l l}
  	[(1+\alpha^2)y -(y_m+d_m\cos\varphi)]\sqrt{x_m^2+y^2} \\ = \alpha[2\cos\varphi y^2 -(d_m+y_m\cos\varphi)y + x_m^2\cos\varphi]
	 \end{array}   \right.  \nonumber 
\end{align}
Taking the square of both sides of the previous equations and combining like terms we obtain the quartic equation shown in \eqref{eq:NonCoopQeq}. \qquad \qquad\qquad \qquad\qquad \qquad\qquad \qquad\qquad \qquad $\square$

\bibliographystyle{IEEEtran}
\bibliography{ReferencesTactics}

\begin{thebibliography}{10}
\providecommand{\url}[1]{#1}
\csname url@samestyle\endcsname
\providecommand{\newblock}{\relax}
\providecommand{\bibinfo}[2]{#2}
\providecommand{\BIBentrySTDinterwordspacing}{\spaceskip=0pt\relax}
\providecommand{\BIBentryALTinterwordstretchfactor}{4}
\providecommand{\BIBentryALTinterwordspacing}{\spaceskip=\fontdimen2\font plus
\BIBentryALTinterwordstretchfactor\fontdimen3\font minus
  \fontdimen4\font\relax}
\providecommand{\BIBforeignlanguage}[2]{{%
\expandafter\ifx\csname l@#1\endcsname\relax
\typeout{** WARNING: IEEEtran.bst: No hyphenation pattern has been}%
\typeout{** loaded for the language `#1'. Using the pattern for}%
\typeout{** the default language instead.}%
\else
\language=\csname l@#1\endcsname
\fi
#2}}
\providecommand{\BIBdecl}{\relax}
\BIBdecl

\bibitem{chandler2004cooperative}
P.~Chandler, ``Cooperative control of a team of uavs for tactical missions,''
  in \emph{AIAA 1st Intelligent Systems Technical Conference}, 2004, p. 6215.

\bibitem{Archibald08}
J.~K. Archibald, J.~C. Hill, N.~A. Jepsen, W.~C. Stirling, and R.~L. Frost, ``A
  satisficing approach to aircraft conflict resolution,'' \emph{IEEE
  Transactions on Systems, Man, and Cybernetics-Part C: Applications and
  Reviews}, vol.~38, no.~4, pp. 510--521, 2008.

\bibitem{wang2017cooperative}
Y.~Wang, E.~Garcia, D.~Casbeer, and F.~Zhang, \emph{Cooperative control of
  multi-agent systems: Theory and applications}.\hskip 1em plus 0.5em minus
  0.4em\relax John Wiley \& Sons, 2017.

\bibitem{Beard02}
R.~W. Beard, T.~W. McLain, M.~Goodrich, and E.~P. Anderson, ``Coordinated
  target assignement and intercept for unmanned air vehicles,'' \emph{IEEE
  Transactions on Robotics and Automation}, vol.~18, no.~6, pp. 911--922, 2002.

\bibitem{Castanon09}
D.~A. Castanon and J.~M. Wohletz, ``Model predictive control for stochastic
  resource allocation,'' \emph{IEEE Transactions on Automatic Control},
  vol.~54, no.~8, pp. 1739--1750, 2009.

\bibitem{Lechevin09}
N.~Lechevin, C.~A. Rabbath, and M.~Lauzon, ``A decision policy for the routing
  and munitions management of multiformations of unmanned combat vehicles in
  adversarial urban environments,'' \emph{IEEE Trans. on Control Systems
  Technology}, vol.~17, no.~3, pp. 505--519, 2009.

\bibitem{Garcia16}
E.~Garcia and D.~W. Casbeer, ``Cooperative task allocation for unmanned
  vehicles with communication delays and conflict resolution,'' \emph{AIAA
  Journal of Aerospace and Information Systems}, vol.~13, no.~2, pp. 1--13,
  2016.

\bibitem{Ernest15}
N.~Ernest, K.~Cohen, E.~Kivelevitch, C.~Shumacher, and D.~Casbeer, ``Genetic
  fuzzy trees and their application towards autonomous training and control of
  a squadron of unmanned combat aerial vehicles,'' \emph{Unmanned Systems},
  vol.~3, no.~3, pp. 185--204, 2015.

\bibitem{galati2008near}
D.~G. Galati and M.~A. Simaan, ``Near-nash targeting strategies for
  heterogeneous teams of autonomous combat vehicles,'' in \emph{Unmanned
  Systems Technology X}, vol. 6962.\hskip 1em plus 0.5em minus 0.4em\relax
  International Society for Optics and Photonics, 2008, p. 69620X.

\bibitem{liu2003application}
Y.~Liu, M.~A. Simaan, and J.~B. Cruz~Jr, ``An application of dynamic nash task
  assignment strategies to multi-team military air operations,''
  \emph{Automatica}, vol.~39, no.~8, pp. 1469--1478, 2003.

\bibitem{galati2007effectiveness}
D.~G. Galati and M.~A. Simaan, ``Effectiveness of the nash strategies in
  competitive multi-team target assignment problems,'' \emph{IEEE Transactions
  on Aerospace and Electronic Systems}, vol.~43, no.~1, pp. 126--134, 2007.

\bibitem{ganapathy2003agreement}
S.~Ganapathy and K.~M. Passino, ``Agreement strategies for cooperative control
  of uninhabited autonomous vehicles,'' in \emph{Proceedings of the 2003
  American Control Conference, 2003.}, vol.~2.\hskip 1em plus 0.5em minus
  0.4em\relax IEEE, 2003, pp. 1026--1031.

\bibitem{Isaacs65}
R.~Isaacs, \emph{Differential Games}.\hskip 1em plus 0.5em minus 0.4em\relax
  New York: Wiley, 1965.

\bibitem{weintraub20}
I.~E. {Weintraub}, M.~{Pachter}, and E.~{Garcia}, ``An introduction to
  pursuit-evasion differential games,'' in \emph{2020 American Control
  Conference (ACC)}, 2020, pp. 1049--1066.

\bibitem{zhou2018efficient}
Z.~Zhou, J.~Ding, H.~Huang, R.~Takei, and C.~Tomlin, ``Efficient path planning
  algorithms in reach-avoid problems,'' \emph{Automatica}, vol.~89, pp. 28--36,
  2018.

\bibitem{chen2014path}
M.~Chen, Z.~Zhou, and C.~J. Tomlin, ``A path defense approach to the
  multiplayer reach-avoid game,'' in \emph{53rd IEEE Conference on Decision and
  Control}.\hskip 1em plus 0.5em minus 0.4em\relax IEEE, 2014, pp. 2420--2426.

\bibitem{Coon17}
M.~Coon and D.~Panagou, ``Control strategies for multiplayer
  target-attacker-defender differential games with double integrator
  dynamics,'' in \emph{56th IEEE Conf. on Decision and Control}, 2017, pp.
  1496--1502.

\bibitem{Bopardikar09}
S.~Bopardikar, F.~Bullo, and J.~P. Hespanha, ``A cooperative homicidal
  chauffeur game,'' \emph{Automatica}, vol.~45, no.~7, pp. 1771--1777, 2009.

\bibitem{Garcia2019}
E.~Garcia, D.~W. Casbeer, and M.~Pachter, ``Design and analysis of
  state-feedback optimal strategies for the differential game of active
  defense,'' \emph{IEEE Transactions on Automatic Control}, vol.~64, no.~2, pp.
  553--568, 2019.

\bibitem{Fuchs17}
Z.~E. Fuchs and P.~P. Khargonekar, ``Generalized engage or retreat differential
  game with escort regions,'' \emph{IEEE Transactions on Automatic Control},
  vol.~62, no.~2, pp. 668--681, 2017.

\bibitem{Pachter17}
M.~Pachter, E.~Garcia, and D.~W. Casbeer, ``Differential game of guarding a
  target,'' \emph{AIAA Journal of Guidance, Control, and Dynamics}, vol.~40,
  no.~11, pp. 2991--2998, 2017.

\bibitem{fisac2015pursuit}
J.~F. Fisac and S.~S. Sastry, ``The pursuit-evasion-defense differential game
  in dynamic constrained environments,'' in \emph{IEEE 54th Annual Conference
  on Decision and Control}, 2015, pp. 4549--4556.

\bibitem{Weintraub20IET}
I.~Weintraub, E.~Garcia, and M.~Pachter, ``Optimal guidance strategy for the
  defense of a non-maneuverable target in 3-dimensions,'' \emph{IET Control
  Theory \& Applications}, vol.~14, no.~11, pp. 1531--1538, 2020.

\bibitem{GarciaCDC18}
E.~Garcia, D.~W. Casbeer, and M.~Pachter, ``The capture-the-flag differential
  game,'' in \emph{IEEE 57th Conference on Decision and Control}, 2018, pp.
  4167--4172.

\bibitem{Garcia21GNC}
E.~Garcia, D.~Tran, D.~W. Casbeer, D.~Milutinovic, and M.~Pachter, ``Beyond
  visual range tactics,'' in \emph{AIAA Guidance, Navigation and Control
  Conference}, 2021.

\bibitem{getz1979qualitative}
W.~M. Getz and G.~Leitmann, ``Qualitative differential games with two
  targets,'' \emph{Journal of Mathematical Analysis and Applications}, vol.~68,
  no.~2, pp. 421--430, 1979.

\bibitem{Getz1981capturability}
W.~M. Getz and M.~Pachter, ``Capturability in a two-target game of two cars,''
  \emph{Journal of Guidance, Control, and Dynamics}, vol.~4, no.~1, pp. 15--21,
  1981.

\bibitem{ardema1985combat}
M.~D. Ardema, M.~Heymann, and N.~Rajan, ``Combat games,'' \emph{Journal of Opt.
  Theory and Applications}, vol.~46, no.~4, pp. 391--398, 1985.

\bibitem{Pachter94}
M.~Pachter, J.~J. D'Azzo, and J.~L. Dargan, ``Automatic formation flight
  control,'' \emph{AIAA Journal of Guidance, Control, and Dynamics}, vol.~17,
  no.~6, pp. 1380--1383, 1994.

\end{thebibliography}

\end{document}